\documentclass[12pt,oneside,reqno]{amsart}
\usepackage{mathrsfs}
\usepackage{graphics}
\usepackage{amssymb}
\usepackage{enumerate}
\newcommand{\dif}{\mathrm{d}}

\newcommand{\be}{\begin{eqnarray}}
\newcommand{\ee}{\end{eqnarray}}
\newcommand{\ce}{\begin{eqnarray*}}
\newcommand{\de}{\end{eqnarray*}}
\newtheorem{theorem}{Theorem}[section]
\newtheorem{lemma}[theorem]{Lemma}
\newtheorem{remark}[theorem]{Remark}
\newtheorem{definition}[theorem]{Definition}
\newtheorem{proposition}[theorem]{Proposition}
\newtheorem{Example}[theorem]{Example}
\newtheorem{corollary}[theorem]{Corollary}
\newtheorem{condition}[theorem]{Condition}
\def\e{\varepsilon}
\def\s{\sigma}
\def\t{\theta}

\def\o{\omega}

\def\d{\delta}
\def\p{\partial}
\def\g{\gamma}

\def\la{\langle}
\def\ra{\rangle}
\def\[{{\Big[}}
\def\]{{\Big]}}
\def\<{{\langle}}
\def\>{{\rangle}}
\def\({{\Big(}}
\def\){{\Big)}}

\def\no{\nonumber}
\def\bt{\begin{theorem}}
\def\et{\end{theorem}}
\def\bl{\begin{lemma}}
\def\el{\end{lemma}}
\def\br{\begin{remark}}
\def\er{\end{remark}}
\def\bx{\begin{Example}}
\def\ex{\end{Example}}
\def\bd{\begin{definition}}
\def\ed{\end{definition}}
\def\bp{\begin{proposition}}
\def\ep{\end{proposition}}
\def\bc{\begin{corollary}}
\def\ec{\end{corollary}}
\def\bco{\begin{condition}}
\def\eco{\end{condition}}
\def\cA{{\mathcal A}}
\def\cB{{\mathcal B}}

\def\cG{{\mathcal G}}

\def\cK{{\mathcal K}}

\def\cM{{\mathcal M}}

\def\cP{{\mathcal P}}

\def\cU{{\mathcal U}}

\def\mE{{\mathbb E}}

\def\mH{{\mathbb H}}

\def\mL{{\mathbb L}}

\def\mN{{\mathbb N}}

\def\mP{{\mathbb P}}
\def\mQ{{\mathbb Q}}
\def\mR{{\mathbb R}}

\def\mU{{\mathbb U}}

\def\mW{{\mathbb W}}

\def\sB{{\mathscr B}}
\def\sC{{\mathscr C}}

\def\sF{{\mathscr F}}

\def\sL{{\mathscr L}}

\def\sS{{\mathscr S}}

\def\geq{\geqslant}
\def\leq{\leqslant}
\def\epsilon{\varepsilon}

\begin{document}

\allowdisplaybreaks

\title{Large deviation principles of nonlinear filtering  for McKean-Vlasov stochastic differential equations}

\author{Huijie Qiao and Shengqing Zhu}

\dedicatory{School of Mathematics,
Southeast University,\\
Nanjing, Jiangsu 211189, P.R.China}



\thanks{This work was partly supported by NSF of China (No.12071071).}

\thanks{Corresponding author: Huijie Qiao, hjqiaogean@seu.edu.cn}


\date{}

\begin{abstract}
In this paper, we study large deviation principles of nonlinear filtering for McKean-Vlasov stochastic differential equations. First of all, we establish the large deviation principle for the space-distribution dependent Zakai equation 
by a weak convergence approach. Then based on the obtained result and the relationship between the space-distribution dependent Zakai equation and the space-distribution dependent Kushner-Stratonovich equation, the large 
deviation principle for the latter is proved.
\end{abstract}

\maketitle \rm

\section{Introduction}

Stochastic filtering deals with the estimation problem under partial information. Given two stochastic processes, the signal process and observation process, the filtering problem aims to estimate the signal process based on the 
observation process. This paper focuses on a general model in which signal and observation processes are McKean-Vlasov  stochastic differential equations (SDEs for short). These equations are different from general SDEs, and 
depend on the positions and probability distributions of particles. In \cite{sc1, sc2} Sen and Caines studied nonlinear filtering problems of McKean-Vlasov SDEs with independent noises. Recently, for McKean-Vlasov SDEs with 
correlated noises, Liu and the first named author \cite{lq} solved the nonlinear filtering problem of them and established the well-posedness of the space-distribution dependent Zakai equation and the space-distribution dependent 
Kushner-Stratonovich equation.

The asymptotic theory of large deviation principles (LDPs for short) aims to make exponential asymptotic estimates of the probability of rare events. LDPs are widely used and have become a very active branch of applied probability. LDPs for nonlinear filtering 
have been established for some cases (cf. \cite{hijab, mpx, mx, pz, xiong2}). Let us mention some works related with our results. In \cite{hijab}, Hijab first studied the LDP for general SDEs. Later, Maroulas and Xiong \cite{mx} generalized the result in \cite{hijab} to SDEs driven by fractional Brownian motions. Recently, Maroulas, Pan and Xiong \cite{mpx} improved the result in \cite{hijab} to SDEs driven by Brownian motions and Poisson random measures. 

For the case in this paper, as far as we know, there are not related LDP results. However, McKean-Vlasov SDEs and related results have been widely applied in many fields, such as biology, game theory and control theory, which motivates us to study the LDP of nonlinear filtering for McKean-Vlasov stochastic differential equations. Here we investigate the LDPs of the space-distribution dependent Zakai equation and the space-distribution dependent Kushner-Stratonovich equation. In order to prove the LDP of the former, the strategy we apply is to prove the Laplace principle, which is equivalent to the large deviation principle, by using a weak convergence argument, proposed in \cite{bdm}. Then based on this result and the relationship between the space-distribution dependent Zakai equation and the space-distribution dependent Kushner-Stratonovich equation, the LDP for the latter is proved.

The rest of this paper is organized as follows. Section \ref{pre} introduces some notations, definitions and a general criterion of large deviations. In Section \ref{twoequ} we give out the space-distribution dependent Zakai equation and the space-distribution dependent Kushner-Stratonovich equation. Section \ref{ldpsta} focuses on the establishment of two LDPs. And the proofs of two LDPs are placed in Section \ref{proozaldp} and \ref{prooksldp}, respectively.

The following convention will be used throughout the paper: $C$, with or without indices, will denote different positive constants whose values may change from one place to another.

\section{Preliminarie}\label{pre}

In this section, we introduce some notations, definitions and a general criterion of large deviations.

\subsection{Notation}\label{nn}
In this subsection, we introduce some notations used in the sequel.
 
Let $|\cdot|$ and $\|\cdot\|$ be norms of vectors and matrices, respectively. Let $A^{*}$ be the transpose of the matrix $A$. 
 
Let $\cB_b(\mR^{n})$ be the set of bounded Borel measurable functions on $\mR^n$. Let $C^{2}(\mR^{n})$ be the space of continuous functions on $\mR^{n}$ which have continuous partial derivatives of order up to $2$, and $C_{b}^{2}(\mR^{n})$ be the collection of functions in $C^2(\mR^n)$ whose derivatives are bounded.

Let $\sB(\mR^{n})$ be the Borel $\s$-field on $\mR^{n}$. Let $\cM(\mR^{n})$ be the set of bounded Borel measures defined on $\sB(\mR^{n})$ carrying the usual topology of weak convergence. Let
$\cP(\mR^{n})$ be the space of probability measures defined on $\sB(\mR^{n})$ and $\cP_{2}(\mR^{n})$ be the collection of the probability measures $\mu$ on $\sB(\mR^{n})$ satisfying
$$
\|\mu\|^{2}:=\int_{\mR^{n}}|x|^{2} \mu(\dif x)<\infty.
$$
We put on $\cP_{2}(\mR^{n})$ a topology induced by the following $2$-Wasserstein metric:
$$
\mW_{2}^{2}(\mu_{1}, \mu_{2}):=\inf _{\pi \in \sC(\mu_{1}, \mu_{2})} \int_{\mR^{n} \times \mR^{n}}|x-y|^{2} \pi(\dif x, \dif y), \quad \mu_{1}, \mu_{2} \in \cP_{2}(\mR^{n}),
$$
where $ \sC(\mu_{1}, \mu_{2})$ denotes the set of the probability measures whose marginal distributions are $\mu_{1}, \mu_{2}$, respectively. It is known that $(\cP_{2}(\mR^{n}), \mW_{2})$ is a Polish space.

\subsection{$L$-derivatives for functions on $\cP_{2}(\mR^{n})$}

In the subsection we recall the definition of $L$-derivative for functions on $\cP_2(\mR^n)$. Let $I$ be the identity map on $\mR^n$. For $\mu\in\cP_2(\mR^n)$ and $\phi\in L^2(\mR^n, \sB(\mR^n), \mu;\mR^n)$, $<\mu,\phi>:=\int_{\mR^n}\phi(x)\mu(\dif x)$. Moreover, by simple calculation, it holds that $\mu\circ(I+\phi)^{-1}\in\cP_2(\mR^n)$.

\bd\label{lderi}
	$(i)$ A function $f: \cP_{2}(\mR^{n}) \mapsto \mR$ is called $L$-differentiable at $\mu \in \cP_{2}(\mR^{n})$, if for $\varphi \in L^{2}(\mR^{n}, \sB(\mR^{n}), \mu ; \mR^{n})$, the functional
	\ce
	\varphi \mapsto f(\mu \circ(I+\varphi)^{-1})
	\de
	is Fr\'echet differentiable at $\varphi=0$; that is, there exists a unique $\g \in L^{2}(\mR^{n}, \sB(\mR^{n}), \mu ; \mR^{n})$ such that \ce\lim \limits_{\la\mu,|\varphi|^{2}\ra \rightarrow 0} \frac{f(\mu \circ(I+\varphi)^{-1})-f(\mu)-<\mu, \g \cdot \varphi>}{\sqrt{\la\mu,|\varphi|^{2}\ra}}=0.\de
	In the case, we denote $\p_{\mu} f(\mu)=\g$ and call it the $L$-derivative of $f$ at $\mu$.
	
	$(ii)$ A function $f: \cP_{2}(\mR^{n}) \mapsto \mR$ is called $L$-differentiable on $\cP_{2}(\mR^{n})$ if $L$-derivative $\p_{\mu} f(\mu)$ exists for all $\mu \in \cP_{2}(\mR^{n})$.
	\ed

\bd
	The function $f$ is said to be in $C^{(1,1)}(\cP_{2}(\mR^{n}))$, if $\p_{\mu} f$ is continuous, for any $\mu \in \cP_{2}(\mR^{n})$, $\p_{\mu} f(\mu)(\cdot)$ is differentiable, and its derivative $\p_{y} \p_{\mu} f: \cP_{2}(\mR^{n}) \times \mR^{n} \mapsto \mR^{n} \otimes \mR^{n}$ is continuous.
\ed

\bd
$(i)$ The function $F: \mR^n\times\cP_2(\mR^n)\mapsto\mR$ is said to be in $C^{2,(1,1)}(\mR^n\times\cP_2(\mR^n))$, if $F(x,\mu)$ is $C^2$ in $x\in\mR^n$ and $C^{(1,1)}$ in $\mu\in\cP_2(\mR^n)$, respectively, and its derivatives 
$$
\partial_x F(x,\mu), \partial^2_x F(x,\mu), \partial_\mu F(x,\mu)(y),  \partial_y\partial_\mu F(x,\mu)(y)
$$ 
are jointly continuous in the corresponding variable family $(x,\mu)$ or $(x,\mu,y)$.
		
$(ii)$ The function $F: \mR^{n} \times \cP_{2}(\mR^{n}) \mapsto \mR$ is said to be in $\sS(\mR^n\times\cP_2(\mR^n))$, if $F \in C^{2,(1,1)}(\mR^{n} \times \cP_{2}(\mR^{n}))$ and for any compact set $\cK \subset \mR^{n} \times \cP_{2}(\mR^{n})$,
$$
\sup\limits_{(x,\mu)\in\cK}\int_{\mR^n}(\|\partial_y\partial_\mu F(x,\mu)(y)\|^2+|\partial_\mu F(x,\mu)(y)|^2)\mu(\dif y)<\infty.
$$

$(iii)$ The function $F: \mR^{n} \times \cP_{2}(\mR^{n}) \mapsto \mR$ is said to be in $C^{2,(1,1)}_c(\mR^n\times\cP_2(\mR^n))$, if $F \in C^{2,(1,1)}(\mR^{n} \times \cP_{2}(\mR^{n}))$ and $F$ has a compact support. 
\ed

It is obvious that $C^{2,(1,1)}_c(\mR^n\times\cP_2(\mR^n))\subset\sS(\mR^n\times\cP_2(\mR^n))\subset C^{2,(1,1)}(\mR^n\times\cP_2(\mR^n))$.

\subsection{A general criterion of large deviations}

The theory of small-noise large deviations concerns with the asymptotic behavior of solutions of stochastic partial differential equations, say $\{X^{\e}\}$, $\e>0$ defined on a probability space $(\Omega, \sF, \mQ)$, which converges exponentially fast as $\e \rightarrow 0$. The decay rate is expressed via a rate function. An equivalent argument of the large deviation principle is the Laplace principle (\cite{drse}). 

\bd
Let $\mU$ be a Polish space. A function $I: \mU \rightarrow[0, \infty]$ is called a rate function on $\mU$, if for each $M<\infty$ the level set $\{\zeta \in$ $\mU: I(\zeta) \leq M\}$ is a compact subset of  $\mU$. 
\ed

\bd 
The family $\{X^{\e}\}$ is said to satisfy the Laplace principle on $\mU$ with the rate function $I$, if for all bounded continuous functions $G$ mapping $\mU$ into $\mR$,
$$
\lim \limits_{\e \rightarrow 0} \e \log \mE\left\{\exp \left[-\frac{1}{\e} G(X^{\e})\right]\right\}=-\inf _{\zeta \in \mU}\left(G(\zeta)+I(\zeta)\right).
$$ 
\ed

Set $\mathbb{H}:=L^{2}([0,T]; \mR^{m})$ and $\|\psi\|_{\mathbb{H}}:=(\int_{0}^{T}|\psi(t)|^{2}\dif t)^{\frac{1}{2}}$ for $\psi\in\mH$. Let $\mathcal{A}$ be the collection of predictable processes $u(\omega, \cdot)$ belonging to $\mathbb{H}$ for $\mQ$-a.s. $\omega$. For each $N\in\mN$ we define two following spaces:
$$
S^{N}:=\left\{\psi \in \mH: \|\psi\|^2_{\mathbb{H}} \leq N\right\}, \quad \cU^{N}= \Big\{u \in \cA: u(\o,\cdot) \in S^{N}, \mQ-a.s.\,\o \Big\}.
$$
We equip $S^{N}$ with the weak convergence topology in $\mH$. So, $S^{N}$ is metrizable as a compact Polish space. In the sequel, $S^{N}$ will be always endowed with this topology.

Next, a set of sufficient conditions for the Laplace principle is presented.

\bco\label{con}
Let $\cG^{\e}: C([0, T], \mathbb{R}^{m}) \rightarrow \mU$ be a family of measurable mappings. There exists a measurable mapping $\cG^{0}: C([0, T], \mathbb{R}^{m}) \rightarrow \mU$ such that the following holds.

$(i)$ For $N \in \mN$ and $\{f_{\e}, \e>0\}\subset S^N, f \in S^{N}$, if $f_{\e} \rightarrow f,$ as $\e \rightarrow 0$, 
$$
\cG^{0}\(\int^{\cdot}_{0} f_{\e}(s) \dif s\) \rightarrow \cG^{0}\(\int^{\cdot}_{0} f(s) \dif s\).
$$

$(ii)$ For $N \in \mN$ and $\{\psi_{\e}, \e>0\}\subset \cU^{N}, \psi \in \cU^{N}$, if $\psi_{\e} \stackrel{d}{\longrightarrow} \psi$, where $\stackrel{d}{\longrightarrow}$ denotes convergence in distribution, as $\e \rightarrow 0$,
$$
\cG^{\e}\(\sqrt{\e} W+\int_{0}^{\cdot} \psi_{\e}(s) \dif s\) \stackrel{d}{\longrightarrow} \cG^{0}\(\int^{\cdot}_{0} \psi(s) \dif s\),
$$
where $W_\cdot$ is a $m$-dimensional Brownian motion.
\eco

For $\zeta \in \mU$, let $S_{\zeta}=\{\psi \in \mH: \zeta=\cG^{0}(\int^{\cdot}_{0} \psi(s) \dif s)\}$. Let $I: \mU \rightarrow[0, \infty]$ be defined by
\ce 
I(\zeta)=\frac{1}{2} \inf _{\psi \in S_{\zeta}}\|\psi\|^2_{\mH}.
\de
By convention, $I(\zeta)=\infty$ if  $S_{\zeta}=\emptyset$.

\bt\label{ldpfun}(\cite[Theorem 4.2]{bdm})
For  $\e>0$, let  $X^{\e}=\cG^{\e}(\sqrt{\e} W)$ and suppose that Condition \ref{con} holds. Then the family $\{X^{\e}\}$ satisfies the large deviation principle on $\mU$ with the rate function $I$ given above.
\et

\section{the space-distribution dependent Zakai equations and the space-distribution dependent Kushner-Stratonovich equations}\label{twoequ}

In this section, we first introduce nonlinear filtering problems for McKean-Vlasov signal-observation systems and present the space-distribution dependent Zakai equations and the space-distribution dependent Kushner-Stratonovich equations. Then we define strong solutions of them and state the well-posedness result about them.

\subsection{Nonlinear filtering problems for McKean-Vlasov signal-observation systems}

Fix $T>0$. Let $(\Omega, \sF,\{\sF_{t}\}_{t \in[0, T]}, \mQ)$ be a complete filtered probability space and $\{B_{t}, t \geq 0 \}$,$\{W_{t}, t \geq 0\}$ be $d$- and $m$-dimensional standard Brownian motions defined on $(\Omega, \sF,\{\sF_{t}\}_{t \in[0, T]}, \mQ)$, respectively. Moreover, $B.$ and $W.$ are mutually independent.
 Consider the following McKean-Vlasov signal-observation system $(X_{t}, Y_{t})$ on $\mR^{n} \times \mR^{m}$ :
\be\left\{\begin{array}{l}
\dif X_{t}=b_{1}(t,X_{t}, \sL_{X_{t}}^\mQ) \dif t+\s_{1}(t,X_{t}, \sL_{X_{t}}^\mQ) \dif B_{t},\quad 0\leq t\leq T,\\
X_0=x_0,\\
\dif Y_{t}=b_{2}(t, X_{t}, \sL_{X_{t}}^\mQ) \dif t+\s_{2}(t) \dif W_{t},\quad 0\leq t\leq T, \\
Y_0=0,
\end{array}
\right.
\label{xysy}
\ee
where $\sL_{X_{t}}^\mQ$ denotes the distribution of $X_{t}$ under the probability measure $\mQ$, and these coefficients $b_{1}: [0,T]\times\mR^{n} \times \cP_{2}(\mR^{n}) \mapsto \mR^{n}$, $\s_{1}: [0,T]\times\mR^{n} \times \cP_{2}(\mR^{n}) \mapsto \mR^{n \times d}$, $b_{2}:[0, T] \times \mR^{n} \times \cP_{2}(\mR^{n}) \mapsto \mR^{m}$, $\s_{2}:[0, T] \mapsto \mR^{m \times m}$ are Borel measurable. 

Assume:
\begin{enumerate}[$(\mathbf{H}_{b_1,\s_1})$]
\item There is an increasing function $L_1: [0,T]\mapsto(0,\infty)$ such that for $t\in[0,T], x_{1}, x_{2} \in \mR^{n},  \mu_{1}, \mu_{2} \in \cP_{2}(\mR^{n})$,
\ce
&&|b_{1}(t,x_{1}, \mu_{1})-b_{1}(t,x_{2}, \mu_{2})|^2+\|\s_{1}(t,x_{1}, \mu_{1})-\s_{1}(t,x_{2}, \mu_{2})\|^{2}\\
&\leq& L_1(t)(|x_{1}-x_{2}|^{2}+\mW_{2}^{2}(\mu_{1}, \mu_{2})).
\de
\end{enumerate}
\begin{enumerate}[$(\mathbf{H}_{b_2,\s_2})$]
\item For any $t \in[0, T]$, $\s_{2}(t)$ is invertible and for any $t \in[0, T], x\in \mR^{n}, \mu \in \cP_{2}(\mR^{n})$
$$
|b_{2}(t, x, \mu)| \vee\|\s_{2}(t)\| \vee\|\s_{2}^{-1}(t)\| \leq L_{2},
$$
where $L_{2}>0$ is a constant.
\end{enumerate}

\br
$(\mathbf{H}_{b_1,\s_1})$ implies that for $t\in[0,T], x \in \mR^{n}, \mu \in \cP_{2}(\mR^{n})$,
$$
|b_{1}(t,x, \mu)|^{2}+\|\s_{1}(t,x, \mu)\|^{2} \leq \bar{L}_{1}(1+|x|+\|\mu\|)^{2},
$$
where $\bar{L}_{1}>0$ is a constant.
\er

Under $(\mathbf{H}_{b_1,\s_1})$ $(\mathbf{H}_{b_2,\s_2})$, Eq.(\ref{xysy}) has a unique strong solution $(X,Y)$. Set 
\ce
\Lambda_{t}^{-1}:=\exp \left\{-\int_{0}^{t} h^{i}(s, X_{s}, \sL_{X_{s}}^\mQ) \dif W^i_{s}-\frac{1}{2} \int_{0}^{t}|h(s, X_{s}, \sL_{X_{s}}^\mQ)|^{2} \dif s\right\},
\de
where $h(t, x, \mu):=\s_{2}^{-1}(t) b_{2}(t, x, \mu)$. Here and hereafter, we use the convention that repeated indices imply summation. Then by $({\bf H}_{b_2,\s_2})$, we know that
$$
\mE\left [\exp \left\{\int_{0}^{T}|h(s, X_{s}, \sL_{X_{s}}^\mQ)|^{2} \dif s\right\}\right ]<\infty,
$$
and $\Lambda_{\cdot}^{-1}$ is an exponential martingale. Define a probability measure $\tilde{\mQ}$ via 
\ce 
\dif \tilde{\mQ}=\Lambda_{T}^{-1} \dif \mQ,
\de 
and under the probability measure $\tilde{\mQ}$,
$$
\tilde{W}_{t}:=W_{t}+\int_{0}^{t} h(s, X_{s}, \sL_{X_{s}}^\mQ) \dif s
$$
is an $(\mathscr{F}_t)$-adapted Brownian motion. Moreover, by $\dif Y_t=\s_{2}(t)\dif \tilde{W}_t$, the $\s$-algebra $\sF_{t}^{Y}$ generated by $\{Y_{s}, 0 \leq s \leq t\}$, can be characterized as
\ce
\sF_{t}^{Y}=\sF_{t}^{\tilde{W}}.
\de
We augment $\mathscr{F}_t^{Y}$ in a usual sense and still denote the augmentation of $\mathscr{F}_t^{Y}$ as $\mathscr{F}_t^Y$. 

Set for any $F\in\sS(\mR^n\times\cP_2(\mR^n))$
\ce
\pi_{t}(F):=\mE[F(X_{t}, \sL_{X_{t}}^\mQ) \mid \sF_{t}^{Y}],
\de 
and the Kallianpur-Striebel formula  gives 
\ce
\pi_{t}(F)=\frac{\tilde{\mP}_{t}(F)}{\tilde{\mP}_{t}(1)},
\de
where 
\ce
\tilde{\mP}_{t}(F):=\tilde{\mE}[F(X_{t}, \sL_{X_{t}}^\mQ) \Lambda_{t} \mid \sF_{t}^{Y}], 
\de
and $\tilde{\mE}$ denotes the expectation under the probability measure $\tilde{\mQ}$. Now we are ready to present the space-distribution dependent Zakai equations and space-distribution dependent Kushner-Stratonovich equations (\cite{lq}).

\bl\label{zakailem} (The space-distribution dependent Zakai equation)
Under $({\bf H}_{b_1,\s_1})$ $({\bf H}_{b_2,\s_2})$, for $F\in\sS(\mR^n\times\cP_2(\mR^n))$, $\tilde{\mP}_{t}(F)$ satisfies the following equation
\be
\tilde{\mP}_{t}(F)=F(x_0,\d_{x_0})+\int_{0}^{t} \tilde{\mP}_{s}(\mL_s F) \dif s+\int_{0}^{t} \tilde{\mP}_{s}[F h^{i}(s, \cdot, \cdot)] \dif \tilde{W}_{s}^{i},
\label{zakai}
\ee
where the operator $\mL_s$ is defined as:
\ce
(\mL_s F)(x,\mu)&=&\p_{x_{i}} F(x, \mu) b_{1}^{i}(s,x, \mu)+\frac{1}{2} \p_{x_{i} x_{j}}^{2} F(x, \mu)(\s_{1} \s_{1}^{*})^{i j}(s,x, \mu)\no\\
&&+\int_{\mR^{n}}(\p_{\mu} F)_{i}(x, \mu)(u) b_{1}^{i}(s,u, \mu) \mu(\dif u)\no\\
&&+\frac{1}{2} \int_{\mR^{n}} \p_{u_{i}}(\p_{\mu} F)_{j}(x, \mu)(u)(\s_{1} \s_{1}^{*})^{i j}(s,u, \mu) \mu(\dif u).
\de
\el

\bl\label{kslem} (The space-distribution dependent Kushner-Stratonovich equation)
Under $({\bf H}_{b_1,\s_1})$ $({\bf H}_{b_2,\s_2})$, for $F\in\sS(\mR^n\times\cP_2(\mR^n))$, $\pi_{t}(F)$ satisfies the following equation
\be
\pi_{t}(F)=F(x_0,\d_{x_0})+\int_{0}^{t} \pi_{s}(\mL_s F) \dif s+\int_{0}^{t}\Big\{\pi_{s}[F h^{i}(s, \cdot, \cdot)]-\pi_{s}(F) \pi_{s}[h^{i}(s, \cdot, \cdot)]\Big\} \dif \hat{W}_{s}^{i},
\label{ks}
\ee
where $\hat{W}_{t}:=\tilde{W}_{t}-\int_{0}^{t} \pi_{s}[h(s, \cdot, \cdot)] \dif s$ is an innovation process.
\el

\subsection{The well-posedness for strong solutions of Eq.(\ref{zakai}) and Eq.(\ref{ks})}

\bd\label{zakastrosolu}
A strong solution for the space-distribution dependent Zakai equation (\ref{zakai}) is a $(\mathscr{F}_t^Y)_{t\in[0,T]}$-adapted, continuous and $\cM(\mR^n\times\cP_2(\mR^n))$-valued process $(\Sigma_t)_{t\in[0,T]}$ such that $(\Sigma_t)_{t\in[0,T]}$
solves the space-distribution dependent Zakai equation (\ref{zakai}), that is,
\be
&&\Sigma_t(F)=F(x_0,\d_{x_0})+\int_0^t\Sigma_s(\mL_s F)\dif s
+\int_0^t\Sigma_s(Fh^j(s,\cdot,\cdot))\dif \tilde{W}^j_s\no\\
&&\qquad\qquad\qquad F\in\sS(\mR^n\times\cP_2(\mR^n)),\quad t\in[0,T].
\label{zakastrosolueq}
\ee
\ed

\bd\label{ksstrosolu}
A strong solution for the space-distribution dependent Kushner-Stratonovich equation (\ref{ks}) is a $(\mathscr{F}_t^Y)_{t\in[0,T]}$-adapted, continuous and $\cP(\mR^n\times\cP_2(\mR^n))$-valued process $(\Pi_t)_{t\in[0,T]}$ such that $(\Pi_t)_{t\in[0,T]}$ solves the space-distribution dependent Kushner-Stratonovich equation (\ref{ks}), that is,
\be
&&\Pi_t(F)=F(x_0,\d_{x_0})+\int_0^t\Pi_s(\mL_s F)\dif s+\int_0^t\(\Pi_s(Fh^j(s,\cdot,\cdot))-\Pi_s(F)\Pi_s(h^j(s,\cdot,\cdot))\)\dif \hat{W}^{\prime j}_s, \no\\
&&\qquad\qquad\qquad\qquad F\in\sS(\mR^n\times\cP_2(\mR^n)), \qquad t\in[0,T], 
\label{kseqstrsolueq}
\ee
where $\hat{W}^{\prime}_t:=\tilde{W}_t-\int_0^t\Pi_s(h(s,\cdot,\cdot))\dif s$.
\ed

Here we state two results about the well-posedness of Eq.(\ref{zakai}) and Eq.(\ref{ks}) (cf. \cite{lq}).

\bt\label{zaun}
Suppose that $({\bf H}_{b_1,\s_1})$ $({\bf H}_{b_2,\s_2})$ hold. Then $\tilde{\mP}_{\cdot}$ is the unique strong solution of the space-distribution dependent Zakai equation (\ref{zakai}). 
\et

\bt\label{ksun}
Suppose that $({\bf H}_{b_1,\s_1})$ $({\bf H}_{b_2,\s_2})$ hold. Then $\pi_{\cdot}$ is the unique strong solution of the space-distribution dependent Kushner-Stratonovich equation (\ref{ks}). 
\et

\section{The Large deviation principles}\label{ldpsta}

In this section we study the limiting behavior of the optimal filter with a small signal-to-noise ratio, i.e., consider the signal given in the framework and the observation process below, for $0<\e<1$,
\ce
Y_{t}^{\e}=\sqrt{\e} \int_{0}^{t} b_{2}(s, X_{s}, \sL_{X_{s}}^\mQ) \dif s+\int_{0}^{t} \s_{2}(s) \dif W_{s}.
\de
Set 
$$
(\Lambda_{t}^{\e})^{-1}:=\exp \left\{-\sqrt{\e} \int_{0}^{t} h^{i}(s, X_{s}, \sL_{X_{s}}^\mQ) \dif W^i_{s}-\frac{\e}{2} \int_{0}^{t}|h(s, X_{s}, \sL_{X_{s}}^\mQ)|^{2} \dif s\right\},
$$
and we define the probability measure $\tilde{\mQ}^{\e}$ by
 \ce
 \mathrm{d} \tilde{\mQ}^{\e}=(\Lambda_{T}^{\e})^{-1} \dif\mQ.
 \de 
 Moreover, we define for any $F\in\sS(\mR^n\times\cP_2(\mR^n))$
 \ce
&&\tilde{\mP}_{t}^{\e}(F):=\tilde{\mE}^{\e}[F(X_{t}, \sL_{X_{t}}^\mQ) \Lambda_{t}^{\e} \mid \sF_{t}^{Y^{\e}}],\\
&&\pi_t^\e(F):=\mE[F(X_{t}, \sL_{X_{t}}^\mQ)\mid \sF_{t}^{Y^{\e}}],
\de 
where $\tilde{\mE}^{\e}$ denotes the expectation under the measure $\tilde{\mQ}^{\e}$. Thus, by the same deduction to that of Eq.(\ref{zakai}) and Eq.(\ref{ks}), one can obtain the following proposition.

\bp 
Under $({\bf H}_{b_1,\s_1})$ $({\bf H}_{b_2,\s_2})$,  for any $F\in\sS(\mR^n\times\cP_2(\mR^n))$, $\tilde{\mP}_{t}^{\e}(F)$ satisfies the following equation
\be
\tilde{\mP}_{t}^{\e}(F)=F(x_0,\d_{x_0})+\int_{0}^{t} \tilde{\mP}_{s}^{\e}(\mL_s F) \dif s+\sqrt{\e}\int_{0}^{t} \tilde{\mP}_{s}^{\e}[F h^{i}(s, \cdot, \cdot)] \dif \tilde{W}_{s}^{\e, i},
\label{ezakai}
\ee
where 
\ce
\tilde{W}_{t}^{\e}:={W}_{t}+\sqrt{\e}\int_{0}^{t} h(s, X_{s}, \sL_{X_{s}}^\mQ)\dif s.
\de
And $\pi_t^\e(F)$ solves the following equation
\be
\pi_{t}^{\e}(F)=F(x_0,\d_{x_0})+\int_{0}^{t} \pi_{s}^{\e}(\mL_s F) \dif s+\sqrt{\e}\int_{0}^{t}\Big\{\pi_{s}^{\e}[F h^{i}(s, \cdot, \cdot)]-\pi_{s}^{\e}(F) \pi_{s}^{\e}[h^{i}(s, \cdot, \cdot)]\Big\} \dif \hat{W}_{s}^{\e, i},
\label{eks}
\ee
where $\hat{W}_{t}^{\e}:=\tilde{W}_{t}^{\e}-\sqrt{\e} \int_{0}^{t} \pi_{s}^{\e}[h(s, \cdot, \cdot)]\dif s$ is an innovation process.
\ep

Now, it is the position to state two main results in this paper.

\bt\label{zaldp} 
Under $({\bf H}_{b_1,\s_1})$ $({\bf H}_{b_2,\s_2})$, $\{\tilde{\mP}^{\e}\}$ satisfies the LDP on $C([0, T], \cM(\mR^{n}\times\cP_2(\mR^n)))$ with the rate function $I_{1}$ given by
\ce
 I_{1}(\zeta):= \frac{1}{2}\inf\limits_{\{\psi \in S_{\zeta}: \zeta=\tilde{\mP}^{0, \psi}\}}\|\psi\|^2_{\mH}, \quad \zeta\in C([0, T], \cM(\mR^{n}\times\cP_2(\mR^n))),
\de
where $\tilde{\mP}^{0, \psi}$ solves the following equation
\ce
\tilde{\mP}_{t}^{0, \psi}(F)=F(x_0,\d_{x_0})+\int_{0}^{t} \tilde{\mP}_{s}^{0, \psi}(\mL_s F) \dif s+\int_{0}^{t} \tilde{\mP}_{s}^{0, \psi}[F h^i(s, \cdot, \cdot)] \psi^i(s) \dif s,~F\in C^{2,(1,1)}_c(\mR^n\times\cP_2(\mR^n)).
\de
\et 

\bt \label{ksldp}
Under $({\bf H}_{b_1,\s_1})$ $({\bf H}_{b_2,\s_2})$, $\{\pi^{\e}\}$ satisfies the LDP on $C([0, T], \cP(\mR^{n}\times\cP_2(\mR^n)))$ with the rate function ${I}_{2}$ given by
\ce
I_{2}(\xi):=\frac{1}{2}\inf\limits_{\{\psi \in S_{\xi}: \xi=\pi^{0, \psi}\}}\|\psi\|^2_{\mH}, \quad \xi\in C([0, T], \cP(\mR^{n}\times\cP_2(\mR^n))),
\de
where $\pi^{0, \psi}$ is the solution of the following equation
\ce
\pi_{t}^{0, \psi}(F)&=&F(x_0,\d_{x_0})+\int_{0}^{t} \pi_{s}^{0, \psi}(\mL_s F)\dif s+\int_{0}^{t}[\pi_{s}^{0, \psi}(Fh^i(s, \cdot, \cdot))\no\\
&&-\pi_{s}^{0, \psi}(F)\pi_{s}^{0, \psi}(h^i(s, \cdot, \cdot))] \psi^i(s)\dif s,~F\in C^{2,(1,1)}_c(\mR^n\times\cP_2(\mR^n)).
\de 
\et

The proofs of two above theorems are placed in Section \ref{proozaldp} and \ref{prooksldp}, respectively.

\section{Proof of Theorem \ref{zaldp}}\label{proozaldp}

In this section, we prove Theorem \ref{zaldp}.

By Theorem \ref{zaun}, we know that $\tilde{\mP}_{\cdot}^{\e}$ is the unique strong solution of Eq.(\ref{ezakai}). Thus, there is a measurable mapping $\cG^{\e}: C([0, T], \mR^{m}) \rightarrow C([0, T], \cM(\mR^{n}\times\cP_2(\mR^n)))$ such that
\ce
\tilde{\mP}_{\cdot}^{\e}:=\cG^{\e}(\sqrt{\e} \tilde{W}^{\e}).
\de
In order to prove the Laplace principle for $\tilde{\mP}_{\cdot}^{\e}$, we will verify Condition \ref{con} with $\mathbb{U}=C([0,T],\cM(\mR^{n}\times\cP_2(\mR^n)))$.

Let $\psi \in \cU^{\mathrm{N}}$. The controlled version of Eq.(\ref{ezakai}) is given by
\be
\tilde{\mP}_{t}^{\e, \psi}(F)&=&F(x_0,\d_{x_0})+\int_{0}^{t} \tilde{\mP}_{s}^{\e, \psi}(\mL_s F) \dif s+\sqrt{\e}\int_{0}^{t} \tilde{\mP}_{s}^{\e, \psi}[F h^{i}(s, \cdot, \cdot)] \dif \tilde{W}_{s}^{\e, i} \no\\
&&+\int_{0}^{t} \tilde{\mP}_{s}^{\e, \psi}[F h^i(s, \cdot, \cdot)] \psi^i(s) \dif s, \quad F\in C^{2,(1,1)}_c(\mR^n\times\cP_2(\mR^n)).
\label{psezakai} 
\ee
The following lemma establishes the existence and uniqueness of Eq.(\ref{psezakai}).

\bl
Suppose that $({\bf H}_{b_1,\s_1})$ $({\bf H}_{b_2,\s_2})$ hold. For $\e>0$, set 
\ce 
\tilde{\mP}^{\e, \psi}:=\cG^{\e}\(\sqrt{\e} \tilde{W}^{\e}+\int^{\cdot}_{0} \psi(s) \dif s\).
\de
Then $\tilde{\mP}^{\e, \psi}$ is the unique strong solution of Eq.(\ref{psezakai}).
\el
\begin{proof} 
Take a control $\psi \in \cU^{N}$, and consider
$$
M_{t}^{\e, \psi}:=\exp \left\{-\frac{1}{\sqrt{\e}} \int_{0}^{t} \psi(s) \dif \tilde{W}_{s}^{\e}-\frac{1}{2 \e} \int_{0}^{t}|\psi(s)|^{2} \dif s\right\}.
$$
Since $\psi \in \cU^{N}$, the Novikov condition holds, and $M_{\cdot}^{\e, \psi}$ is an exponential martingale. Define a new probability measure $\tilde{\mQ}^{\e, \psi}$ by 
\ce 
\dif \tilde{\mQ}^{\e, \psi}=M_{T}^{\e, \psi} \dif \tilde{\mQ}^{\e},
\de 
and $\tilde{\mQ}^{\e, \psi}$ is an equivalent probability measure with respect to $\tilde{\mQ}^{\e}$. By Girsanov's theorem,
\ce 
\tilde{W}_{t}^{\e, \psi}:=\tilde{W}_{t}^{\e}+\frac{1}{\sqrt{\e}} \int_{0}^{t} \psi(s)\dif s 
\de
is a Brownian motion with respect to $\tilde{\mQ}^{\e, \psi}$. Hence, the desired result follows from replacing $\tilde{W}$ in Eq.(\ref{zakai}) with $\tilde{W}^{\e, \psi}$ and Theorem \ref{zaun}.
\end{proof}

Now, consider the controlled version of Eq.(\ref{ezakai}) without noises, i.e. for $F\in C^{2,(1,1)}_c(\mR^n\times\cP_2(\mR^n)), \psi\in S^{N}$
\be
\tilde{\mP}_{t}^{0, \psi}(F)=F(x_0,\d_{x_0})+\int_{0}^{t} \tilde{\mP}_{s}^{0, \psi}\(\mL_s F+h^i(s, \cdot, \cdot)\psi^i(s)F\)\dif s.
\label{zakainono}
\ee

\bl\label{ps0zakaiwell}
Suppose that $({\bf H}_{b_1,\s_1})$ $({\bf H}_{b_2,\s_2})$ hold. Then Eq.(\ref{zakainono}) has a unique solution $\tilde{\mP}^{0, \psi}$.
\el
\begin{proof}
{\bf Existence.} Set
\ce
\Lambda_{t}^{0, \psi}:=\exp \left\{\int_{0}^{t} h^i(s, X_{s}, \sL_{X_{s}}^\mQ) \psi^i(s) \dif s\right\},\quad \tilde{\mP}_{t}^{0, \psi}(F):=\mE[\Lambda_{t}^{0, \psi} F(X_{t}, \sL_{X_{t}}^\mQ)],
\de
and $\tilde{\mP}^{0, \psi}$ is a solution of Eq.(\ref{zakainono}). Indeed, on one hand, it holds that
$$
\dif \Lambda_{t}^{0, \psi}=\Lambda_{t }^{0, \psi} h^i(t,X_{t}, \sL_{X_{t}}^\mQ) \psi^i(t) \dif t.
$$
On the other hand, by the It\^o formula, for any $F\in C^{2,(1,1)}_c(\mR^n\times\cP_2(\mR^n))$, we have that
\ce
F(X_{t}, \sL_{X_{t}}^\mQ)=F(x_0,\d_{x_0})+\int_{0}^{t} \mL_s F(X_{s}, \sL_{X_{s}}^\mQ) \dif s+ \int_{0}^{t} \p_{x_{i}} F(X_{s}, \sL_{X s}^{\mQ}) \s_{1}^{i j}(X_{s}, \sL_{X s}^{\mQ}) \dif W_{s}^{j}.
\de
Then the It\^o formula together with the above deduction yields that
\ce
\Lambda_{t}^{0, \psi} F(X_{t}, \sL_{X_{t}}^\mQ)&=&F(x_0,\d_{x_0})+\int_{0}^{t} \Lambda_{s}^{0, \psi}F(X_{s}, \sL_{X s}^{\mQ})h^i(s, X_{s}, \sL_{X_{s}}^\mQ) \psi^i(s)\dif s\\
&&+\int_{0}^{t} \Lambda_{s}^{0, \psi}\mL_s F(X_{s}, \sL_{X_{s}}^\mQ)\dif s+\int_{0}^{t} \Lambda_{s}^{0, \psi}\p_{x_{i}} F(X_{s}, \sL_{X s}^{\mQ}) \s_{1}^{i j}(X_{s}, \sL_{X s}^{\mQ}) \dif W_{s}^{j}.
\de
Taking the expectation on both sides of the above equality, we obtain that $\tilde{\mP}_{t}^{0, \psi}$ is a solution of Eq.(\ref{zakainono}).

{\bf Uniqueness.} Theorem 3.5 (A slightly extended version) in \cite{bk1}  implies the uniqueness of solutions for Eq.(\ref{zakainono}).
\end{proof}

Define the measurable mapping
\ce
\cG^{0}: C([0, T], \mR^{m}) \rightarrow C([0, T], \cM(\mR^{n}\times\cP_2(\mR^n)))
\de 
by
\ce
\cG^{0}\(\int_{0}^{\cdot} \psi(s) \dif s\)=\tilde{\mP}^{0, \psi}.
\de
Then we verify Condition \ref{con} through $\cG^{\e}, \cG^{0}$.

\bl\label{vericond2}
Suppose that $\psi_{\e}, \psi\in \cU^{N}$, $\psi_{\e}\rightarrow\psi$ almost surely as $\e \rightarrow 0$. Then it holds that
$$
\cG^{\e}\(\sqrt{\e} \tilde{W}^{\e}+\int^{\cdot}_{0} \psi_{\e}(s) \dif s\)\stackrel{\tilde{\mQ}^{\e}}{\rightarrow} \cG^{0}\(\int^{\cdot}_{0} \psi(s) \dif s\).
$$
\el
\begin{proof} 
{\bf Step 1.} We prove that the family
$$
\left\{\tilde{\mP}^{\e, \psi_{\e}}=\cG^{\e}\(\sqrt{\e} \tilde{W}^{\e}+\int^{\cdot}_{0} \psi_{\e}(s) \dif s\), \e \in(0,1)\right\}
$$
is tight in $C([0, T], \cM(\mR^{n}\times\cP_2(\mR^n)))$.

It is well-known, e.g. see \cite{kx}  that we only need to prove the tightness of $\{\tilde{\mP}^{\e, \psi_{\e}}(F); \e\in (0,1)\}$ in $C([0, T], \mathbb{R})$ for every test function $F\in C^{2,(1,1)}_c(\mR^n\times\cP_2(\mR^n))$. In fact, we only need to justify that 

$(i)$ $\sup\limits_{\e\in (0,1)}\tilde{\mE}^\e\sup\limits _{t \in[0, T]} \left|\tilde{\mP}_{t}^{\e, \psi_{\e}}(F)\right|\leq C$;

$(ii)$ For any $(\mathscr{F}_t)_{t\geq0}$-stopping time $\tau\leq T$
\ce
\limsup_{\t\downarrow0}\sup\limits_{\e\in(0,1)}\sup\limits_{0\leq\tau<\tau+\t\leq T}\tilde{\mE}^\e|\tilde{\mP}_{\tau+\t}^{\e, \psi_{\e}}(F)-\tilde{\mP}_{\tau}^{\e, \psi_{\e}}(F)|=0.
\de

Note that $\tilde{\mP}_{t}^{\e, \psi_{\e}}(1)$ satisfies the following equation
\ce
\tilde{\mP}_{t}^{\e, \psi_{\e}}(1)=1+\sqrt{\e}\int_{0}^{t} \tilde{\mP}_{s}^{\e, \psi_{\e}}[h^{i}(s, \cdot, \cdot)] \dif \tilde{W}_{s}^{\e, i}+\int_{0}^{t} \tilde{\mP}_{s}^{\e, \psi_{\e}}[h^i(s, \cdot, \cdot)] \psi^i_{\e}(s) \dif s.
\de
Thus, by It\^o's formula, we get that
\ce
\(\tilde{\mP}_{t}^{\e, \psi_{\e}}(1)\)^{2}&=&1+2 \sqrt{\e} \int_{0}^{t} \tilde{\mP}_{s}^{\e, \psi_{\e}}(1) \tilde{\mP}_{s}^{\e, \psi_{\e}}[h^{i}(s, \cdot, \cdot)] \dif \tilde{W}_{s}^{\e, i}\\
&&+2 \int_{0}^{t} \tilde{\mP}_{s}^{\e, \psi_{\e}}(1) \tilde{\mP}_{s}^{\e, \psi_{\e}}[h^i(s, \cdot, \cdot)] \psi^i_{\e}(s) \dif s+\e \sum_{i=1}^m\int_{0}^{t}|\tilde{\mP}_{s}^{\e, \psi_{\e}}[h^i(s, \cdot, \cdot)]|^{2} \dif s \\
&=:&1+A_{t}^{1}+A_{t}^{2}+A_{t}^{3}.
\de

For $A_{t}^{1}$, by the BDG inequality, it holds that
\ce
\tilde{\mE}^{\e}\sup\limits_{t \in[0, T]}|A_{t}^{1}|&\leq& 2C\sqrt{\e}\tilde{\mE}^{\e}\left(\sum_{i=1}^m\int_{0}^{T} \(\tilde{\mP}_{s}^{\e, \psi_{\e}}(1)\)^2\left|\tilde{\mP}_{s}^{\e, \psi_{\e}}[h^i(s, \cdot, \cdot)]\right|^2 \dif s\right)^{1/2}\\
&\leq&\frac{1}{4}\tilde{\mE}^{\e}\sup\limits_{t \in[0, T]}\(\tilde{\mP}_{t}^{\e, \psi_{\e}}(1)\)^{2}+C\sum_{i=1}^m\int_{0}^{T}\tilde{\mE}^{\e}\left|\tilde{\mP}_{s}^{\e, \psi_{\e}}[h^i(s, \cdot, \cdot)]\right|^2 \dif s\\
&\leq&\frac{1}{4}\tilde{\mE}^{\e}\sup\limits_{t \in[0, T]}\(\tilde{\mP}_{t}^{\e, \psi_{\e}}(1)\)^{2}+C\int_{0}^{T}\tilde{\mE}^{\e}\(\tilde{\mP}_{s}^{\e, \psi_{\e}}(1)\)^{2} \dif s,
\de
where $(\mathbf{H}_{b_2,\s_2})$ is used in the last inequality.

For $A_{t}^{2}$, the H\"older inequality implies that
\ce 
\tilde{\mE}^{\e}\sup\limits_{t \in[0, T]}|A_{t}^{2}|&\leq&2\tilde{\mE}^{\e}\sup\limits_{t \in[0, T]}\tilde{\mP}_{t}^{\e, \psi_{\e}}(1)\left(\sum_{i=1}^m\int_{0}^{T} |\tilde{\mP}_{s}^{\e, \psi_{\e}}[h^i(s, \cdot, \cdot)]||\psi^i_{\e}(s)|\dif s\right)\\
&\leq&\frac{1}{4}\tilde{\mE}^{\e}\sup\limits_{t \in[0, T]}\(\tilde{\mP}_{t}^{\e, \psi_{\e}}(1)\)^{2}+C\tilde{\mE}^{\e}\left(\sum_{i=1}^m\int_{0}^{T} |\tilde{\mP}_{s}^{\e, \psi_{\e}}[h^i(s, \cdot, \cdot)]||\psi^i_{\e}(s)|\dif s\right)^2 \no\\
&\leq&\frac{1}{4}\tilde{\mE}^{\e}\sup\limits_{t \in[0, T]}\(\tilde{\mP}_{t}^{\e, \psi_{\e}}(1)\)^{2}+C\tilde{\mE}^{\e}\sum_{i=1}^m\int_{0}^{T} |\tilde{\mP}_{s}^{\e, \psi_{\e}}[h^i(s, \cdot, \cdot)]|^2\dif s \int_{0}^{T}|\psi^i_{\e}(s)|^2\dif s\\ 
&\leq&\frac{1}{4}\tilde{\mE}^{\e}\sup\limits_{t \in[0, T]}\(\tilde{\mP}_{t}^{\e, \psi_{\e}}(1)\)^{2}+C\int_{0}^{T}\tilde{\mE}^{\e}\(\tilde{\mP}_{s}^{\e, \psi_{\e}}(1)\)^{2} \dif s.
\de

For $A_{t}^{3}$, from $(\mathbf{H}_{b_2,\s_2})$, it follows that
\ce 
\tilde{\mE}^{\e}\sup\limits_{t \in[0, T]}|A_{t}^{3}|\leq C\int_{0}^{T}\tilde{\mE}^{\e}\(\tilde{\mP}_{s}^{\e, \psi_{\e}}(1)\)^{2} \dif s.
\de

Combining the above deduction, we obtain that
$$
\tilde{\mE}^{\e}\sup\limits_{t \in[0, T]}\(\tilde{\mP}_{t}^{\e, \psi_{\e}}(1)\)^{2} \leq 2+C\int_{0}^{T}\tilde{\mE}^{\e}\(\tilde{\mP}_{s}^{\e, \psi_{\e}}(1)\)^{2} \dif s.
$$
Then the Gronwall inequality implies that
\be
\sup _{\e \in(0,1)} \tilde{\mE}^{\e}\sup\limits_{t \in[0, T]}\(\tilde{\mP}_{t}^{\e, \psi_{\e}}(1)\)^{2} \leq C,
\label{pepsi1}
\ee
which together with the boundness of $F$ yields that
\ce
\sup _{\e \in(0,1)} \tilde{\mE}^{\e}\(\sup _{t \in[0, T]} \(\tilde{\mP}_{t}^{\e, \psi_{\e}}(F)\)^{2}\) \leq C.
\de
Thus, we complete the proof of $(i)$.

Next, we verify $(ii)$. Note that $\tilde{\mP}_{t}^{\e, \psi_{\e}}(F)$ satisfies Eq.(\ref{psezakai}), i.e.
\ce
\tilde{\mP}_{t}^{\e, \psi_\e}(F)&=&F(x_0,\d_{x_0})+\int_{0}^{t} \tilde{\mP}_{s}^{\e, \psi_\e}(\mL_s F) \dif s+\sqrt{\e}\int_{0}^{t} \tilde{\mP}_{s}^{\e, \psi_\e}[F h^{i}(s, \cdot, \cdot)] \dif \tilde{W}_{s}^{\e, i} \no\\
&&+\int_{0}^{t} \tilde{\mP}_{s}^{\e, \psi_\e}[F h^i(s, \cdot, \cdot)] \psi_\e^i(s) \dif s.
\de
Thus, for any stopping time $\tau$ and any $\t>0$, by the H\"older inequality and the BDG inequality, it holds that
\ce
\tilde{\mE}^\e|\tilde{\mP}_{\tau+\t}^{\e, \psi_{\e}}(F)-\tilde{\mP}_{\tau}^{\e, \psi_{\e}}(F)|^2&\leq& 3\tilde{\mE}^\e\left|\int_{\tau}^{\tau+\t} \tilde{\mP}_{s}^{\e, \psi_\e}(\mL_s F) \dif s\right|^2+3\e\tilde{\mE}^\e\left|\int_{\tau}^{\tau+\t} \tilde{\mP}_{s}^{\e, \psi_\e}[F h^{i}(s, \cdot, \cdot)] \dif \tilde{W}_{s}^{\e, i}\right|^2\\
&&+3\tilde{\mE}^\e\left|\sum_{i=1}^m\int_{\tau}^{\tau+\t} \tilde{\mP}_{s}^{\e, \psi_\e}[F h^i(s, \cdot, \cdot)] \psi_\e^i(s) \dif s\right|^2\\
&\leq&3\t\tilde{\mE}^\e\int_{\tau}^{\tau+\t} |\tilde{\mP}_{s}^{\e, \psi_\e}(\mL_s F)|^2 \dif s+3\e\tilde{\mE}^\e\sum_{i=1}^m\int_{\tau}^{\tau+\t} |\tilde{\mP}_{s}^{\e, \psi_\e}[F h^{i}(s, \cdot, \cdot)] |^2\dif s\\
&&+3m\tilde{\mE}^\e\sum_{i=1}^m\int_{\tau}^{\tau+\t}\left|\tilde{\mP}_{s}^{\e, \psi_{\e}}[F h^i(s, \cdot, \cdot)]\right|^{2} \dif s\int_{\tau}^{\tau+\t}|\psi^i_{\e}(s)|^2 \dif s\\
&\leq&3\t^{2} \tilde{\mE}^{\e}\(\underset{t \in[0, T]}{\sup }\(\tilde{\mP}_{t}^{\e, \psi_{\e}}(\mL_s F)\)^{2}\)+3\e\tilde{\mE}^\e\sum_{i=1}^m\int_{\tau}^{\tau+\t} |\tilde{\mP}_{s}^{\e, \psi_\e}[F h^{i}(s, \cdot, \cdot)] |^2\dif s\\
&&+3C\tilde{\mE}^\e\sum_{i=1}^m\int_{\tau}^{\tau+\t}\left|\tilde{\mP}_{s}^{\e, \psi_{\e}}[F h^i(s, \cdot, \cdot)]\right|^{2} \dif s\\
&\leq&3\t^{2} \tilde{\mE}^{\e}\(\underset{t \in[0, T]}{\sup }\(\tilde{\mP}_{t}^{\e, \psi_{\e}}(|\mL_s F|)\)^{2}\)+C\t\tilde{\mE}^\e\(\underset{t \in[0, T]}{\sup }\(\tilde{\mP}_{t}^{\e, \psi_{\e}}(|F||h(t, \cdot, \cdot)|)\)^{2}\)\\
&\leq&C(\t^{2}+\t).
\de
From this, we obtain that
\ce
\limsup_{\t\downarrow0}\sup\limits_{\e\in(0,1)}\sup\limits_{0\leq\tau<\tau+\t\leq T}\tilde{\mE}^\e|\tilde{\mP}_{\tau+\t}^{\e, \psi_{\e}}(F)-\tilde{\mP}_{\tau}^{\e, \psi_{\e}}(F)|=0,
\de
which is just $(ii)$.

{\bf Step 2.} We show that a weak limit of $\{\tilde{\mP}^{\e, \psi_{\e}}, \e>0\}$ is $\tilde{\mP}^{0, \psi}=\cG^{0}(\int_{0}^{\cdot} \psi(s) \dif s)$.

By {\bf Step 1.}, there exists a $\bar{\tilde{\mP}}^{0, \psi}\in C([0, T], \cM(\mR^{n}\times\cP_2(\mR^n)))$ such that $\tilde{\mP}^{\e, \psi_{\e}}$ converges weakly to $\bar{\tilde{\mP}}^{0, \psi}$. Therefore, we only need to prove that $\bar{\tilde{\mP}}^{0, \psi}=\tilde{\mP}^{0, \psi}$.

First of all, $\tilde{\mP}_{t}^{\e, \psi_{\e}}(F)$ solves the following equation
\ce
\tilde{\mP}_{t}^{\e, \psi_{\e}}(F)&=&\tilde{\mP}_{0}^{\e, \psi_{\e}}(F)+\int_{0}^{t} \tilde{\mP}_{s}^{\e, \psi_{\e}}(\mL_s F) \dif s+\sqrt{\e}\int_{0}^{t} \tilde{\mP}_{s}^{\e, \psi_\e}[F h^{i}(s, \cdot, \cdot)] \dif \tilde{W}_{s}^{\e, i}\\
&&+\int_{0}^{t} \tilde{\mP}_{s}^{\e, \psi_{\e}}[F h^i(s, \cdot, \cdot)] \psi^i_{\e}(s) \dif s.
\de
On one side, by the BDG inequality, it holds that
\ce
\tilde{\mE}^{\e}\left(\sup _{t \in[0, T]}\left|\int_{0}^{t} \tilde{\mP}_{s}^{\e, \psi_\e}[F h^{i}(s, \cdot, \cdot)] \dif \tilde{W}_{s}^{\e, i}\right|^{2}\right) &\leq& C\tilde{\mE}^{\e} \sum_{i=1}^m\int_{0}^{T} |\tilde{\mP}_{s}^{\e, \psi_\e}[F h^{i}(s, \cdot, \cdot)]|^2 \dif s\\
&\leq&CT\tilde{\mE}^\e \sum_{i=1}^m\(\underset{t \in[0, T]}{\sup }\(\tilde{\mP}_{t}^{\e, \psi_{\e}}(F h^{i}(s, \cdot, \cdot))\)^{2}\)\\
&\leq&CTm\tilde{\mE}^\e\underset{t \in[0, T]}{\sup }\(\tilde{\mP}_{t}^{\e, \psi_{\e}}(1)\)^2.
\de
Therefore, the dominated convergence implies that as $\e$ tends to $0$,
\ce
\sqrt{\e}\int_{0}^{t} \tilde{\mP}_{s}^{\e, \psi_\e}[F h^{i}(s, \cdot, \cdot)] \dif \tilde{W}_{s}^{\e, i}\rightarrow 0, \quad \tilde{\mQ}^{\e}-a.s..
\de

On the other side, we note that
\ce
\int_{0}^{t} |\tilde{\mP}_{s}^{\e, \psi_{\e}}[F h^i(s, \cdot, \cdot)]|^2\dif s\leq C\int_{0}^{t} |\tilde{\mP}_{s}^{\e, \psi_{\e}}(1)|^2\dif s\leq CT\underset{t \in[0, T]}{\sup }\(\tilde{\mP}_{t}^{\e, \psi_{\e}}(1)\)^2<\infty, \tilde{\mQ}^{\e}-a.s..
\de
Thus, by the assumption that $\psi_{\e}\rightarrow\psi$ almost surely as $\e \rightarrow 0$, it holds that as $\e \rightarrow 0$,
\ce
\int_{0}^{t} \tilde{\mP}_{s}^{\e, \psi_{\e}}[F h^i(s, \cdot, \cdot)] \psi^i_{\e}(s) \dif s\rightarrow \int_{0}^{t} \bar{\tilde{\mP}}_{s}^{0, \psi}[F h^i(s, \cdot, \cdot)] \psi^i(s) \dif s, \tilde{\mQ}^{\e}-a.s..
\de

Finally, as $\e \rightarrow 0$, since $\tilde{\mP}_{t}^{\e, \psi_{\e}}(F)\rightarrow\bar{\tilde{\mP}}_{t}^{0, \psi}(F)$, and $\tilde{\mP}_{0}^{\e, \psi_{\e}}(F)\rightarrow\bar{\tilde{\mP}}_{0}^{0, \psi}(F)$, $\bar{\tilde{\mP}}^{0, \psi}$ solves Eq.(\ref{zakainono}). Then the uniqueness of the solution to Eq.(\ref{zakainono}) completes the proof.
\end{proof}

By the similar or even simpler deduction to that of Lemma \ref{vericond2}, we obtain the following result.

\bl\label{vericond1} 
Suppose that $f_{\e}, f\in S^{N}$, $f_{\e} \rightarrow f$, as $\e \rightarrow 0$. Then it holds that
$$
\cG^{0}\(\int^{\cdot}_{0} f_{\e} (s) \dif s\) \rightarrow \cG^{0}\(\int^{\cdot}_{0} f(s) \dif s\).
$$ 
\el

Now, it is the position to prove Theorem \ref{zaldp}.

{\bf Proof of Theorem \ref{zaldp}.}

By Theorem \ref{ldpfun}, to establish LDP, it is sufficient to verify the two conditions in Condition \ref{con}.  
In Lemma \ref{vericond1}, we have already proved Condition \ref{con} $(i)$. Thus we only need to verify Condition \ref{con} $(ii)$. 

For $\epsilon\in(0,1)$ and $\{\psi_{\epsilon}, \e>0\}\subset\cU^{N}$, $\psi\in\cU^{N}$, let $\psi_{\epsilon}$ converge to $\psi$ in distribution. By the Skorohod theorem, there exists a
probability space $(\hat{\Omega}, \hat{\sF}, \hat{\mQ})$, and $S^N$-valued random variables $\{\hat{\psi}_{\epsilon}\}$, $\hat{\psi}$ and a $m$-dimensional Brownian motion $\hat{\tilde{W}}$ such that

(i) $\sL_{(\hat{\psi}_{\epsilon},\hat{\tilde{W}})}=\sL_{(\psi_{\epsilon},\tilde{W}^\e)}$ and $\sL_{\hat{\psi}}=\sL_{\psi}$;

(ii) $\hat{\psi}_{\epsilon}$ converges to $\hat{\psi}$ almost surely.

In the following, we construct two Zakai equations:
\be
\hat{\tilde{\mP}}_{t}^{\e, \hat{\psi}_\e}(F)&=&F(x_0,\d_{x_0})+\int_{0}^{t} \hat{\tilde{\mP}}_{s}^{\e, \hat{\psi}_\e}(\mL_s F) \dif s+\sqrt{\e}\int_{0}^{t} \hat{\tilde{\mP}}_{s}^{\e, \hat{\psi}_\e}[F h^{i}(s, \cdot, \cdot)] \dif \hat{\tilde{W}}_{s}^{i} \no\\
&&+\int_{0}^{t} \hat{\tilde{\mP}}_{s}^{\e, \hat{\psi}_\e}[F h^i(s, \cdot, \cdot)] \hat{\psi}_\e^i(s) \dif s,\label{contanal1}\\
\hat{\tilde{\mP}}_{t}^{0, \hat{\psi}}(F)&=&F(x_0,\d_{x_0})+\int_{0}^{t} \hat{\tilde{\mP}}_{s}^{0, \hat{\psi}}(\mL_s F) \dif s+\int_{0}^{t} \hat{\tilde{\mP}}_{s}^{0, \hat{\psi}}[F h^i(s, \cdot, \cdot)] \hat{\psi}^i(s) \dif s.
 \label{deteequa1}
\ee
For Eq.(\ref{contanal1}) and Eq.(\ref{deteequa1}), each has a unique strong solution denoted by $\hat{\tilde{\mP}}^{\e, \hat{\psi}_\e}$ and $\hat{\tilde{\mP}}^{0, \hat{\psi}}$ respectively. Moreover, it holds that
$$
\hat{\tilde{\mP}}^{\e, \hat{\psi}_\e}=\cG^{\epsilon}(\sqrt{\epsilon}\hat{\tilde{W}}+\int_{0}^{\cdot}\hat{\psi}_{\epsilon}(s)\dif s), \quad \hat{\tilde{\mP}}^{0, \hat{\psi}}=\cG^{0}(\int_{0}^{\cdot}\hat{\psi}(s)\dif s).
$$
By Lemma \ref{vericond2}, we have that $\cG^{\e}\(\sqrt{\e} \hat{\tilde{W}}+\int^{\cdot}_{0} \hat{\psi}_{\e}(s) \dif s\)\stackrel{\tilde{\mQ}^{\e}}{\rightarrow} \cG^{0}\(\int^{\cdot}_{0} \hat{\psi}(s) \dif s\)$, which yields that
$\cG^{\e}\(\sqrt{\e}  \hat{\tilde{W}}+\int^{\cdot}_{0} \hat{\psi}_{\e}(s) \dif s\)\rightarrow \cG^{0}\(\int^{\cdot}_{0} \hat{\psi}(s) \dif s\)$ in distribution. Note that
\ce
\cG^{\e}\(\sqrt{\e} \tilde{W}^{\e}+\int^{\cdot}_{0} \psi_{\e}(s) \dif s\)&\overset{d}{=}&\cG^{\epsilon}(\sqrt{\epsilon}\hat{\tilde{W}}+\int_{0}^{\cdot}\hat{\psi}_{\epsilon}(s)\dif s),\\
\cG^{0}\(\int^{\cdot}_{0} \psi(s) \dif s\)&\overset{d}{=}&\cG^{0}\(\int^{\cdot}_{0} \hat{\psi}(s) \dif s\).
\de
So, $\cG^{\e}\(\sqrt{\e} \tilde{W}^{\e}+\int^{\cdot}_{0} \psi_{\e}(s) \dif s\)\rightarrow \cG^{0}\(\int^{\cdot}_{0} \psi(s) \dif s\)$ in distribution, which is Condition \ref{con} $(ii)$.

\section{Proof of Theorem \ref{ksldp}}\label{prooksldp}

In this section, we prove Theorem \ref{ksldp}. First of all, we consider the control version of Eq.(\ref{eks}) without noises, i.e. for $F\in C^{2,(1,1)}_c(\mR^n\times\cP_2(\mR^n)), \psi\in S^N$
\be
\pi_{t}^{0, \psi}(F)&=&F(x_0,\d_{x_0})+\int_{0}^{t} \pi_{s}^{0, \psi}(\mL_s F)\dif s\no\\
&&+\int_{0}^{t}[\pi_{s}^{0, \psi}(Fh^i(s, \cdot, \cdot))-\pi_{s}^{0, \psi}(F)\pi_{s}^{0, \psi}(h^i(s, \cdot, \cdot))] \psi^i(s)\dif s.
\label{ksnono}
\ee

\bl\label{zakaiks}
Suppose that $({\bf H}_{b_1,\s_1})$ $({\bf H}_{b_2,\s_2})$ hold. Then Eq.(\ref{ksnono}) has a unique solution $\pi^{0, \psi}$.
\el
\begin{proof}
{\bf Existence.} We claim that $\pi_{t}^{0, \psi}(F):=\frac{\tilde{\mP}_{t}^{0, \psi}(F)}{\tilde{\mP}_{t}^{0, \psi}(1)}$ satisfies Eq.(\ref{ksnono}), where $\tilde{\mP}^{0, \psi}(F)$ solves Eq.(\ref{zakainono}). Indeed, it holds that
\ce 
\tilde{\mP}_{t}^{0, \psi}(F)=F(x_0,\d_{x_0})+\int_{0}^{t} \tilde{\mP}_{s}^{0, \psi}(\mL_s F) \dif s+\int_{0}^{t} \tilde{\mP}_{s}^{0, \psi}[F h^i(s, \cdot, \cdot)] \psi^i(s) \dif s,
\de
and furthermore
\ce 
\tilde{\mP}_{t}^{0, \psi}(1)=1+\int_{0}^{t} \tilde{\mP}_{s}^{0, \psi}[h^i(s, \cdot, \cdot)] \psi^i(s) \dif s,
\de
which implies that
\ce 
\dif \frac{1}{\tilde{\mP}_{t}^{0, \psi}(1)}=-\frac{1}{\(\tilde{\mP}_{t}^{0, \psi}(1)\)^{2}}\tilde{\mP}_{t}^{0, \psi}[h^i(t, \cdot, \cdot)] \psi^i(t) \dif t.
\de
So, by the Taylor formula we get that
\ce 
\frac{\tilde{\mP}_{t}^{0, \psi}(F)}{\tilde{\mP}_{t}^{0, \psi}(1)}
&=&F(x_0,\d_{x_0})+\int_{0}^{t} \frac{\tilde{\mP}_{s}^{0, \psi}(\mL_s F)}{\tilde{\mP}_{s}^{0, \psi}(1)} \dif s+\int_{0}^{t} \frac{\tilde{\mP}_{s}^{0, \psi}[F h^i(s, \cdot, \cdot)] }{\tilde{\mP}_{s}^{0, \psi}(1)}\psi^i(s)\dif s \\
&&-\int_{0}^{t} \frac{\tilde{\mP}_{s}^{0, \psi}[h^i(s, \cdot, \cdot)]}{\tilde{\mP}_{s}^{0, \psi}(1)}\frac{\tilde{\mP}_{s}^{0, \psi}(F)}{\tilde{\mP}_{s}^{0, \psi}(1)}\psi^i(s)\dif s.
\de
That is, $\frac{\tilde{\mP}_{t}^{0, \psi}(F)}{\tilde{\mP}_{t}^{0, \psi}(1)}$ satisfies Eq.(\ref{ksnono}).

{\bf Uniqueness.} Assume that $\pi^{0, \psi,1}, \pi^{0, \psi,2}$ are two solutions of Eq.(\ref{ksnono}). So, by Lemma \ref{kszakai} below, it holds that $\tilde{\mP}_{t}^{0, \psi,j}(F):=\pi_{t}^{0, \psi,j}(F) M_{t}^{\psi,j}$ is 
the solution for Eq.(\ref{zakainono}) for $j=1,2$ and 
$$
\pi_{t}^{0, \psi,j}(F)=\frac{\tilde{\mP}_{t}^{0, \psi,j}(F)}{\tilde{\mP}_{t}^{0, \psi,j}(1)}.
$$
Besides, by Lemma \ref{ps0zakaiwell}, we know that $\tilde{\mP}_{t}^{0, \psi,1}(F)=\tilde{\mP}_{t}^{0, \psi,2}(F)$ for $F\in C^{2,(1,1)}_c(\mR^n\times\cP_2(\mR^n))$. Hence, $\pi_t^{0, \psi,1}(F)=\pi_{t}^{0, \psi,2}(F)$ for 
$F\in C^{2,(1,1)}_c(\mR^n\times\cP_2(\mR^n))$. The proof is complete.
\end{proof}

\bl\label{kszakai} 
Assume that $\pi^{0, \psi}$ is a solution of Eq.(\ref{ksnono}). Then it holds that $\pi^{0, \psi}(F) M^{\psi}$ is a solution of Eq.(\ref{zakainono}), where $M_{t}^{\psi}:=\exp \left\{\int_{0}^{t} \pi^{0, \psi}_{s}[h^i(s, \cdot, \cdot)] \psi^i(s) \dif s\right\}$. 
\el
\begin{proof}
First, we know that
\ce
&&\dif M_{t}^{\psi}=M_{t}^{\psi}\pi^{0, \psi}_{t}[h^i(t, \cdot, \cdot)]\psi^i(t)\dif t, \\
&&\dif \pi_{t}^{0, \psi}(F)=\pi_{t}^{0, \psi}(\mL_s F)\dif t+[\pi_{t}^{0, \psi}(Fh^i(t, \cdot, \cdot))-\pi_{t}^{0, \psi}(F)\pi_{t}^{0, \psi}(h^i(t, \cdot, \cdot))] \psi^i(t)\dif t.
\de
Thus, the Taylor formula implies that
\ce
\dif M_{t}^{\psi} \pi_{t}^{0, \psi}(F)&=&\pi_{t}^{0, \psi}(F) M_{t}^{\psi} \pi^{0, \psi}_{t}[h^i(t, \cdot, \cdot)]\psi^i(t)\dif t+\pi_{t}^{0, \psi}(\mL_s F) M_{t}^{\psi} \dif t \no\\
&&+\pi^{0, \psi}_{t}[F h^i(t, \cdot, \cdot)] \psi^i(t) M_{t}^{\psi} \dif t-\pi_{t}^{0, \psi}(F) M_{t}^{\psi} \pi^{0, \psi}_{t}[h^i(t, \cdot, \cdot)] \psi^i(t) \dif t\no\\
&=&\pi_{t}^{0, \psi}(\mL_s F) M_{t}^{\psi} \dif t+\pi^{0, \psi}_{t}[F h^i(t, \cdot, \cdot)] \psi^i(t) M_{t}^{\psi} \dif t,
\de
and furthermore
\ce
M_{t}^{\psi} \pi_{t}^{0, \psi}(F)=F(x_0,\d_{x_0})+\int_{0}^{t}\pi_{s}^{0, \psi}(\mL_s F) M_{s}^{\psi} \dif s+\int_{0}^{t}\pi^{0, \psi}_{s}[F h^i(s, \cdot, \cdot)]M_{s}^{\psi}\psi^i(s) \dif s,
\de
which coincides  with Eq.(\ref{zakainono}). Therefore, the uniqueness of Eq.(\ref{zakainono}) yields that
\ce
\tilde{\mP}_{t}^{0, \psi}(F)=\pi_{t}^{0, \psi}(F) M_{t}^{\psi}.
\de
The proof is complete.
\end{proof}

Now, we are ready to prove Theorem \ref{ksldp}.

{\bf Proof of Theorem \ref{ksldp}.}

First of all, we define a mapping $\Sigma: C([0, T],\cM(\mR^{n}\times \cP_{2}(\mR^{n})) \setminus\{0\}) \rightarrow C([0, T], \cP(\mR^{n}\times \cP_{2}(\mR^{n})))$ as follows
\ce
\Sigma(\tilde{\mP})_t(\cdot):=\frac{\tilde{\mP}_{t}(\cdot)}{\tilde{\mP}_{t}(1)}, \quad \tilde{\mP}\in C([0, T],\cM(\mR^{n}\times \cP_{2}(\mR^{n})) \setminus\{0\}).
\de
Then by Theorem \ref{zaldp} and the contraction principle, it holds that $\pi^{\e}=\Sigma(\tilde{\mP}^\e)$ satisfies the large deviation principle with the rate function given by
\be
\hat{I}_2(\xi)=\inf\limits_{\xi=\Sigma(\zeta)} I_1(\zeta), \quad \xi\in C([0, T], \cP(\mR^{n}\times \cP_{2}(\mR^{n}))).
\label{hati2}
\ee
In order to obtain the required result, we only need to show that $I_{2}=\hat{I}_{2}$.

If $I_2(\xi)<\infty$, for any $\d>0$ there exists a $\psi\in S_\xi$ such that $\xi=\pi^{0,\psi}$ and 
$$
\frac{1}{2}\int_{0}^{T} \psi^{2}(s) \dif s<I_2(\xi)+\d.
$$
Besides, by Lemma \ref{kszakai}, it holds that $\tilde{\mP}_{t}^{0, \psi}=\pi_{t}^{0, \psi} M_{t}^{\psi}$. Thus, Lemma \ref{zakaiks} implies that $\xi=\pi^{0,\psi}=\Sigma(\tilde{\mP}^{0, \psi})$, which together with (\ref{hati2}) yields that
\be
\hat{I}_2(\xi)\leq I_1(\tilde{\mP}^{0, \psi})\leq \frac{1}{2}\int_{0}^{T} \psi^{2}(s) \dif s<I_2(\xi)+\d.
\label{dayu}
\ee

Next, if $\hat{I}_2(\xi)<\infty$, for any $\d>0$ there exists a $\zeta\in C([0, T],\cM(\mR^{n}\times \cP_{2}(\mR^{n})))$ such that $\xi=\Sigma(\zeta)$ and $I_1(\zeta)\leq \hat{I}_2(\xi)+\d$. Take $\psi\in S_{\zeta}$ such that $\tilde{\mP}^{0,\psi}=\zeta$ and 
\ce
\frac{1}{2} \int_{0}^{T} \psi^{2}(s) \dif s<I_{1}(\zeta)+\d.
\de
Besides, by Lemma \ref{zakaiks}, $\xi=\Sigma(\zeta)=\Sigma(\tilde{\mP}^{0,\psi})=\pi^{0,\psi}$. Therefore, $\psi\in S_{\xi}$ and 
\be
I_2(\xi)\leq\frac{1}{2}\int_{0}^{T} \psi^{2}(s) \dif s<I_{1}(\zeta)+\d<\hat{I}_2(\xi)+2\d.
\label{xiyu}
\ee

Finally, combining (\ref{xiyu}) and (\ref{dayu}) and letting $\d\rightarrow 0$, we obtain that $I_{2}=\hat{I}_{2}$.

\end{document}